\documentclass[oneside]{amsart}
\newcommand{\formatswitch}{preprint}

\usepackage{amssymb}
\usepackage{ifthen}
\usepackage{verbatim}
\usepackage{psfig}
\usepackage[all]{xy}

\newcommand{\tref}[1]{(\ref{#1})}

\ifthenelse{\equal{\formatswitch}{big}}{
\setlength{\textwidth}{432pt}
\setlength{\oddsidemargin}{18.8775pt}

\setcounter{tocdepth}{1}
}{
\ifthenelse{\equal{\formatswitch}{paper}}{

\setcounter{tocdepth}{1}
}{
\setcounter{tocdepth}{2}
}
}
\DeclareMathAlphabet\EuScript{U}{eus}{m}{n}
\DeclareMathAlphabet\EuScriptb{U}{eus}{b}{n}

\newcommand{\sscr}[1]{\EuScript{#1}}

\newcommand{\claimenum}{\renewcommand{\theenumi}{\alph{enumi}}
 \renewcommand{\labelenumi}{\textit{(\theenumi)}}
 \renewcommand{\theenumii}{\roman{enumii}}
 \renewcommand{\labelenumii}{\textit{(\theenumii)}}
 \begin{enumerate}}
\newcommand{\claimenumend}{\end{enumerate}}

\newcommand{\romanenum}{\renewcommand{\theenumi}{\roman{enumi}}
 \renewcommand{\labelenumi}{\textit{(\theenumi)}}
 \renewcommand{\theenumii}{\alph{enumii}}
 \renewcommand{\labelenumii}{\textit{(\theenumii)}}
 \begin{enumerate}}
\newcommand{\romanenumend}{\end{enumerate}}

\newtheorem{dummy}{realdumb}[section]
\newtheorem{thm}{Theorem}

\newtheorem{lemma}[dummy]{Lemma}
\newtheorem{prop}[dummy]{Proposition}
{\theoremstyle{definition} }

{\theoremstyle{definition} }

\newtheorem{cora}[dummy]{Corollary}

\renewcommand{\text}{\mathrm}

\newcommand{\strutdepth}{\dp\strutbox}
\newcommand{\marginalnote}[1]
   {\strut\vadjust{\kern-\strutdepth\domarginalnote{#1}}}
\newcommand{\domarginalnote}[1]{\vtop to \strutdepth{
  \baselineskip\strutdepth
   \vss\llap{ #1\ \ }\null}}  
\newcounter{showlabelflag}
\newcounter{makelabelflag}
\newcommand{\showlabels}{\setcounter{showlabelflag}{1}}

\newcommand{\makelabels}{\setcounter{makelabelflag}{1}}
\newcommand{\hidelabels}{\setcounter{showlabelflag}{2}}
\newcommand{\mylabel}[1]{
  \ifthenelse{\value{makelabelflag}=1}
    {\label{#1}}{}
  \ifthenelse{\value{showlabelflag}=1}
    {\marginpar{#1}}{}\relax}

\newcommand{\N}{{\mathbf N}}
\newcommand{\sub}{\subseteq}

\ifthenelse{\equal{\formatswitch}{draft}}{
\showlabels
}{\relax}

\newcommand{\ZS}{Zappa-Sz\'ep }

\newcommand{\scr}{\sscr}

\newcommand{\zapprod}{\mathbin{\bowtie}}
\newcommand{\FB}{\scr{F}\zapprod B_\infty}
\newcommand{\FP}{\scr{F}\zapprod S_\infty}

\newcommand{\mymargin}[1]{
  \ifthenelse{\value{showlabelflag}=1}
    {\marginpar{#1}}{}\relax}

\newcounter{enumo}

\newcommand{\RRsh}{\kern -1 pt \Rsh}

\newcounter{keepitemnum}

\newcounter{keepitemnumm}

\begin{document}

\bibliographystyle{amsplain}
\begin{center}{\bfseries The Algebra of Strand Splitting. II.
\\
A Presentation for the Braid Group on One Strand\footnote{AMS
Classification (2000): primary 20F05, secondary 20F36, 20E99, 
20B07}}\end{center}
\vspace{3pt}
\begin{center}{MATTHEW G. BRIN}\end{center}
\vspace{4pt}
\vspace{3pt}
\begin{center}May 28, 2004\end{center}

\CompileMatrices


\makelabels
\hidelabels


\section{Introduction}\mylabel{IntroSec}

In \cite{brin:bv}, we give descriptions of a braided version \(BV\)
of Thompson's group \(V\) as well as a group \(\widehat{BV}\) that
contains \(BV\) as a subgroup and that is somewhat easier to work
with.  The paper \cite{brin:bv} contains both geometric and
algebraic descriptions of these two groups and shows that for each
group the two descriptions are of isomorphic groups.  An infinite
presentation for \(\widehat{BV}\) is also given in \cite{brin:bv}.
The current paper computes a finite presentation for
\(\widehat{BV}\) as well as two infinite and two finite
presentations for \(BV\).

There is a very close relationship (similar to that between an Artin
group and its corresponding Coxeter group) between \(BV\) and \(V\)
as well as between \(\widehat{BV}\) and a corresponding group
\(\widehat{V}\).  The relationship is close enough that, with no
extra work, our calculations also yield finite and infinite
presentations of \(V\) and finite presentations of
\(\widehat{V}\).  Infinite presentations of \(\widehat{V}\) are
given in \cite{brin:bv}.

The similarity of the relationships to that of an Artin group to its
corresponding Coxeter group is expressed by the fact that
presentations for \(V\) and \(\widehat{V}\) can be obtained from the
presentations for \(BV\) and \(\widehat{BV}\) by simply declaring
the squares of certain generators to be the identity.  We give
multiple presentations of the groups \(V\) and \(BV\) since one set
more closely resembles existing presentations of Thompson's groups
(see, for example, \cite{CFP}), while the other set emphasizes the
Artin/Coxeter relation since the presentations (finite and infinite)
for \(V\) are gotten from the corresponding presentations for \(BV\)
by declaring the squares of \emph{all} the generators to be the
identity.  See Corollary \ref{InfiniteBVPresC} and the last
paragraph of Theorem \ref{FiniteBVPres}.

The finite presentation for \(V\) given in \cite{CFP} has fewer
relations than the finite presentation that we obtain for \(V\), but
the presentation in \cite{CFP} is harder to relate to the
presentation of \(BV\).

We work entirely with the algebraic structure of the groups above,
and we refer the reader to \cite{brin:bv} for the geometric
descriptions that reveal the structures as braid groups.  The title
refers to the fact that \(\widehat{BV}\) can be viewed as a braid
group on countably many strands that are allowed to split and
recombine as long as the order of splits and recombinations is
remembered.  The subgroup \(BV\) corresponds to the subgroup in
which all splitting and recombining is confined to the first strand
and in which all braiding is confined to the strands obtained from
splitting the first strand.  Thus \(BV\) can be thought of as the
``braid group on one strand'' where splitting and joining is
allowed.

In the next section, we list what we need from \cite{brin:bv} and
also give needed facts about another of Thompson's groups.

\section{Introducing the groups}

\subsection{The groups \protect\(\widehat{V}\protect\) and
\protect\(\widehat{BV}\protect\)}

The following presentations come from Theorem 1 of \cite{brin:bv}
where \(\Lambda = \{\lambda_0, \lambda_1, \ldots\}\) and
\(\Sigma=\{\sigma_0, \sigma_1, \ldots\}\).
\mymargin{InfinitePresI--II}
\begin{alignat}{2} \label{InfinitePresI}
\widehat{V} = 
  \langle
    \Lambda \cup \Sigma \mid \,
    & \lambda_q\lambda_m = \lambda_m \lambda_{q+1}, 
      &\qquad&m<q, \\ \notag
    &\sigma_m^2=1, 
      &&m\ge0, \\ \notag
    &\sigma_m\sigma_n = \sigma_n\sigma_m, 
      &&|m-n|\ge2, \\ \notag
    &\sigma_m\sigma_{m+1}\sigma_m = \sigma_{m+1}\sigma_m\sigma_{m+1},
      &&m\ge0, \\ \notag
    &\sigma_q\lambda_m = \lambda_m\sigma_{q+1},
      &&m<q, \\ \notag
    &\sigma_m\lambda_m = \lambda_{m+1}\sigma_m\sigma_{m+1},
      &&m\ge0, \\ \notag
    &\sigma_m\lambda_{m+1} = \lambda_m\sigma_{m+1}\sigma_m,
      &&m\ge0, \\ \notag
    &\sigma_q\lambda_m = \lambda_m\sigma_q, 
      &&m>q+1
  \rangle, \\ \notag \\ \label{InfinitePresII}
\widehat{BV} =
  \langle 
    \Lambda \cup \Sigma \mid \, 
    &\lambda_q\lambda_m = \lambda_m \lambda_{q+1}, 
      &&m<q, \\ \notag
    &\sigma_m\sigma_n = \sigma_n\sigma_m,
      &&|m-n|\ge2, \\ \notag
    &\sigma_m\sigma_{m+1}\sigma_m = \sigma_{m+1}\sigma_m\sigma_{m+1},
      &&m\ge0, \\ \notag
    &\sigma^\epsilon_q\lambda_m = \lambda_m\sigma^\epsilon_{q+1},
      &&m<q,\,\,\epsilon=\pm1, \\ \notag
    &\sigma^\epsilon_m\lambda_m = 
          \lambda_{m+1}\sigma^\epsilon_m\sigma^\epsilon_{m+1},
      &&m\ge0,\,\,\epsilon=\pm1, \\ \notag
    &\sigma^\epsilon_m\lambda_{m+1} = 
          \lambda_m\sigma^\epsilon_{m+1}\sigma^\epsilon_m,
      &&m\ge0,\,\,\epsilon=\pm1, \\ \notag
    &\sigma^\epsilon_q\lambda_m = \lambda_m\sigma^\epsilon_q, 
      &&m>q+1,\,\,\epsilon=\pm1
  \rangle. 
\end{alignat} 
Some of the relations are redundant.  The relations
\(\sigma_m\lambda_{m+1} = \lambda_m\sigma_{m+1}\sigma_m\) follow
from the relations \(\sigma_m\lambda_m =
\lambda_{m+1}\sigma_m\sigma_{m+1}\) in \(\widehat{V}\) by bringing
each \(\sigma_m\) and \(\sigma_{m+1}\) to the other side fo the
equality.  Similarly, the relations \(\sigma^\epsilon_m\lambda_{m+1}
= \lambda_m\sigma^\epsilon_{m+1}\sigma^\epsilon_m\) follow from the
relations \(\sigma^\epsilon_m\lambda_m =
\lambda_{m+1}\sigma^\epsilon_m\sigma^\epsilon_{m+1}\) in
\(\widehat{BV}\).  Also, the exponents \(\epsilon\) can be
eliminated from several of the relations in \(\widehat{BV}\) because
of the group setting.  The relations are listed to be used and not
to give minimal presentations.

The submonoid of \(\widehat{V}\) generated by the positive powers of
the generators in \(\Lambda\) is isomorphic to the submonoid of
\(\widehat{BV}\) generated by the positive powers of generators from
\(\Lambda\) and we will use \(\scr{F}\) to denote both of these
submonoids.  

The submonoid of \(\widehat{V}\) generated by the elements of
\(\Sigma\) and their inverses is isomorphic to the infinite
symmetric group and will be denoted by \(S_\infty\).  The submonoid
of \(\widehat{BV}\) generated by the elements of \(\Sigma\) and
their inverses is isomorphic to the braid group on infinitely many
strands and will be denoted by \(B_\infty\).  We will often have
need to refer separately to the relations of \(B_\infty\), which are
\mymargin{PermRelsB--C} \begin{alignat}{2}\label{PermRelsB}
\sigma_m\sigma_n &= \sigma_n\sigma_m, &\qquad&|m-n|\ge2, \\
\label{PermRelsC}
\sigma_m\sigma_{m+1}\sigma_m&=\sigma_{m+1}\sigma_m\sigma_{m+1},
&&m\ge0. \end{alignat} We regard \(S_\infty\) as the set of finite
permutations of \(\N\), and the generator \(\sigma_i\) corresponds
to the transposition that switches \(i\) and \(i+1\).  The relations
of \(S_\infty\) are \tref{PermRelsB} and \tref{PermRelsC} as well as
\mymargin{PermRelsA}\begin{equation}\label{PermRelsA} \sigma_m^2=1,
\qquad m\ge0. \end{equation} The subgroup of \(S_\infty\) or
\(B_\infty\) generated by \(\{\sigma_0, \sigma_1, \ldots,
\sigma_{n-2}\}\) will be denoted, respectively, by \(S_n\) or
\(B_n\).  That is, \(S_n\) permutes the first \(n\) elements
\(\{0,1,\ldots, n-1\}\) of \(\N\) and \(B_n\) is the braid group on
\(n\) strands.

The submonoid of \(\widehat{V}\) generated by \(\scr{F}\cup
S_\infty\) has the structure of a \ZS product \(\FP\) and the
submonoid of \(\widehat{BV}\) generated by \(\scr{F}\cup B_\infty\)
has the structure of a \ZS product \(\FB\).  The \ZS product is a
generalization of the semidirect product in that neither factor is
required to be normal in the result.

The details of the \ZS product are discussed in \cite{brin:bv} and
in more detail in \cite{brin:zs} and need not concern us here beyond
the following facts.  
{\begin{enumerate}
\renewcommand{\theenumi}{Z\arabic{enumi}} 
\item\mylabel{ZSEmb} The \ZS product \(A\zapprod B\) of two monoids
\(A\) and \(B\) is defined on the set \(A\times B\) and sending all
\(a\in A\) to \((a,1)\) and all \(b\in B\) to \((1, b)\) are
isomorphic embeddings of \(A\) and \(B\) into \(A\zapprod B\).
\item\mylabel{ZSMul} The multiplication on \(A\times B\) gives the
equality \((a,1)(1,b)=(a,b)\) so it makes sense to write \(ab\) for
the pair \((a,b)\) in \(A\times B\).  The reader should note that it
is not always the case that \((1,b)(a,1)=(a,b)\).  
\item\mylabel{ZSUnique} Every element in \(A\zapprod B\) is uniquely
expressible as a product \(ab\) with \(a\in A\) and \(b\in B\).
\end{enumerate}}

The group \(\widehat{V}\) is a group of right fractions of \(\FP\)
and \(\widehat{BV}\) is a group of right fractions of \(\FB\).  What
we need to know about groups of right fractions is the following
where \(M\) is a monoid and \(G\) is its group of right fractions.
See Section 2.3 of \cite{brin:bv} for more details.
{\begin{enumerate}
\renewcommand{\theenumi}{R\arabic{enumi}}
\item\mylabel{GOFEmbed} There is an isomorphic embedding
\(i:M\rightarrow G\). 
\item\mylabel{GOFRep} For every \(g\in G\) there are \(p\) and \(n\)
in \(M\) so that \(g=(i(p))(i(n))^{-1}\).
\end{enumerate}}

\subsection{The monoid \protect\(\scr{F}\protect\) and the group
\protect\(F\protect\)}

The only relations in \(\widehat{V}\) and \(\widehat{BV}\) that
apply to the submonoids called \(\scr{F}\) are those in the first
line of the presentations \tref{InfinitePresI} and
\tref{InfinitePresII}.  These relations have remarkable properties
and we need the following standard facts about these relations as applied to
\(\scr{F}\) and the associated Thompson group \(F\) with
presentation 
\mymargin{GroupFPres}\begin{equation}\label{GroupFPres}
F=\langle \Lambda\mid \lambda_q\lambda_m = \lambda_m \lambda_{q+1},
\quad m<q\rangle 
\end{equation}
The second fact below quotes Proposition 3.3 of \cite{brin:bv} while
the first quotes that proposition and Proposition 6.1 of
\cite{brin:bv}.   The third is 
Lemma 2.1.5 of \cite{brin+fer}
or Theorem 4.3 of \cite{CFP}.
{\begin{enumerate}
\renewcommand{\theenumi}{F\arabic{enumi}}
\item\mylabel{MonoidFPres} The monoid \(\scr{F}\) has presentation given by
\tref{GroupFPres} regarded as a monoid presentation.
\item\mylabel{ForNorm} Each element of \(\scr{F}\) can be written
uniquely as a word 
\(\lambda_{i_0}\lambda_{i_1}\ldots \lambda_{i_k}\) for which
\(i_0\le i_1\le \cdots \le i_k\). 
\item\mylabel{NPNAQ} Every proper quotient of the group \(F\) is
abelian. 
\end{enumerate}}

\subsection{The groups \protect\(V\protect\) and
\protect\(BV\protect\)} 

Section 7 of \cite{brin:bv} identifies \(V\) and \(BV\) as
consisting of those elements of \(\widehat{V}\) and
\(\widehat{BV}\), respectively, of the form
\((F\beta)(G\gamma)^{-1}\) with \(F\) and \(G\) from \(\scr{F}\) and
\(\beta\) and \(\gamma\) from \(S_\infty\) or \(B_\infty\), as
appropriate, that satisfy the following properties for some integer
\(n\ge0\). 
{\begin{enumerate} 
\renewcommand{\theenumi}{S\arabic{enumi}}
\item\mylabel{BalancedCond} The lengths of \(F\) and \(G\) as words
in the generators from \(\Lambda\) are both
equal to \(n\). 
\item\mylabel{SimpleCond} Both \(F\) and \(G\) can be written as
words in the generators 
from \(\Lambda\) in the form \(\lambda_{i_1}\lambda_{i_2} \cdots
\lambda_{i_n}\) so that \(i_j<j\) for all \(j\) with \(1\le j\le
n\).
\item\mylabel{BraidCond} Both \(\beta\) and \(\gamma\) can be
expressed as words in 
elements of \(\{\sigma_0, \sigma_1, \ldots,
\sigma_{n-1}\}\).
\end{enumerate}} 
In \tref{BalancedCond}, the lengths are well defined because of
\tref{MonoidFPres} and the fact that the relations in
\tref{GroupFPres} preserve length.
Note that \tref{BraidCond} says that \(\beta\) and
\(\gamma\) lie in \(S_{n+1}\) or \(B_{n+1}\) as appropriate.

\section{Finite presentations for \protect\(\widehat{V}\protect\)
and \protect\(\widehat{BV}\protect\)} \mylabel{ThompsonArgsSec}

The infinite presentations given in \tref{InfinitePresI} and
\tref{InfinitePresII} fit nicely into a machine due to Thompson that
takes certain structures and reduces them from infinite to finite
\cite{thomp:embed, thompson:cylinder, thompson:adreka}.  We sketch
some details of the machinery by showing how it applies to the
presentations above.  We will use the same machinery to get finite
presentations for \(V\) and \(BV\) after first getting infinite
presentations for those groups.

Because of the group setting, two of the relation sets above can be
written as \(\lambda_m^{-1}\lambda_q\lambda_m = \lambda_{q+1}\) and
\(\lambda_m^{-1}\sigma_q\lambda_m = \sigma_{q+1}\) for all
combinations of \(m\) and \(q\) in which \(m<q\).  Restricting this
to \(m=0\) allows us to write
\(\lambda_i=\lambda_0^{1-i}\lambda_1\lambda_0^{i-1}\) and \(\sigma_i
= \lambda_0^{1-i}\sigma_1\lambda_0^{i-1}\) for \(i>1\).  This shows
that the generating sets can be reduced to \(\{\lambda_0, \lambda_1,
\sigma_0, \sigma_1\}\) for both \(\widehat{V}\) and
\(\widehat{BV}\).  This is not minimal,
but we will ignore that for now.

Using the finite generating sets just given, and using
\(\lambda_i=\lambda_0^{1-i}\lambda_1\lambda_0^{i-1}\) and \(\sigma_i
= \lambda_0^{1-i}\sigma_1\lambda_0^{i-1}\) for \(i>1\) as
definitions, we see that the relations
\(\lambda_q\lambda_m=\lambda_m\lambda_{q+1}\) and
\(\sigma_q\lambda_m=\lambda_m\sigma_{q+1}\) hold by definition for
\(0=m<q\).  We want to know what it takes to conclude the remaining
relations in the presentations in \tref{InfinitePresI} and
\tref{InfinitePresII}.

The relations break into two types: those that involve only one
variable to express all the subscripts involved (such as the ``braid
relation'' \(\sigma_m\sigma_{m+1}\sigma_m =
\sigma_{m+1}\sigma_m\sigma_{m+1}\)) and those that involve two
variables to express all the subscripts involved (such as
\(\sigma_m\sigma_n=\sigma_m\sigma_n\) whenever \(|m-n|\ge2\)).  

A relation of the first type in which the minimum subscript is at
least 1 can be conjugated on both sides by \(\lambda_0^i\) with
\(i>0\) to uniformly elevate all of the subscripts in the relation
by \(i\).  Since all subscripts are non-negative, only two relations
from each such family are needed to generate all the relations in
that family: one in which the subscript 0 appears, and one in which
the minimum subscript is 1.  For example, all of the braid relations
\(\sigma_m\sigma_{m+1}\sigma_m = \sigma_{m+1}\sigma_m\sigma_{m+1}\)
follow from \(\sigma_0\sigma_{1}\sigma_0 =
\sigma_{1}\sigma_0\sigma_{1}\) and \(\sigma_1\sigma_{2}\sigma_1 =
\sigma_{2}\sigma_1\sigma_{2}\).

Relations of the second type can have their appearance altered by
rewriting the larger of the two subscripts.  For example,
\(\sigma_q\lambda_m = \lambda_m\sigma_q\) for \(m>q+1\) can be
written as \(\sigma_q\lambda_{q+k} = \lambda_{q+k}\sigma_q\) for
\(k>1\).  Each value of \(k\) gives a family of relations determined
by the lower subscript.  Again, conjugating by powers of
\(\lambda_0\) allow us to generate all relations for a fixed value
of \(k\) from two relations: one in which the lower subscript is 0,
and one in which the lower subscript is 1.  The task is to get all
values of \(k\) from a finite number of relations.

All of the relations of the second type are statements about the
result of a conjugation: \(X^{-1}YX=Z\).  A commutativity statement
comes under this description as well.  All can be arranged so that the
symbol with the lower subscript can play the role of conjugator,
represented by \(X\).  If the conjugator \(X\) is some \(\sigma_q\)
and \(k\ge3\), then conjugating the relation by \(\lambda^i_{q+2}\)
with \(i>0\) will elevate the value of \(k\) by \(i\) while keeping
the value of \(q\) untouched.  This requires knowing the relations
\(\sigma_q\lambda_{q+2} = \lambda_{q+2}\lambda_q\) for all \(q\) and
knowing inductively the result of conjugating \(\lambda_{i+j}\)
and/or \(\sigma_{i+j}\) by \(\lambda_i\) for values of \(j\)
satisfying \(j<k\).

If the conjugator \(X\) is some \(\lambda_m\) and \(k\ge3\), then
the generator in the role of \(Y\) can be replaced by a conjugate
using a lower subscript.  For example,
\[\lambda_m^{-1}\sigma_{m+k}\lambda_m =
\lambda_m^{-1}(\lambda_{m+1}^{-1}
\sigma_{m+k-1}\lambda_{m+1})\lambda_m.\] Under an inductive
hypothesis based on \(k\) and the homomorphic behavior of
conjugation, the right expression is equal to
\(\lambda_{m+2}^{-1}\sigma_{m+k}\lambda_{m+2}\) which equals
\(\sigma_{m+k+1}\) appealing again to the inductive hypothesis.

Using these arguments and the comments following 
\tref{InfinitePresII} about redundant relations, the reader can verify
that the presentations above can be reduced to the following.  In
the presentations below, we have exploited the fact that one of the
relations we are left with for \(\widehat{V}\) and \(\widehat{BV}\)
is \(\sigma_0\lambda_0 = \lambda_1\sigma_0\sigma_0\).  This allows
us to write \(\lambda_1 = \sigma_0\lambda_0 \sigma_1^{-1}
\sigma_0^{-1}\).

\begin{lemma}\mylabel{FinitePres} The groups 
\(\widehat{V}\) and \(\widehat{BV}\) are presented as groups by the
following in which \(\sigma_i =
\lambda_0^{1-i}\sigma_1\lambda_0^{i-1}\) and
\(\lambda_i=\lambda_0^{1-i}\lambda_1\lambda_0^{i-1}\) for \(i>1\),
as well as  \(\lambda_1 = \sigma_0\lambda_0 \sigma_1^{-1}
\sigma_0^{-1}\)
are definitions: 
\begin{alignat*}{2}
\widehat{V} = 
  \langle
    \lambda_0, \sigma_0, \sigma_1 \mid
      &\lambda_1^{-1}\lambda_2\lambda_1 = \lambda_3,
        &\qquad&\lambda_1^{-1}\lambda_3\lambda_1 = \lambda_4, \\
      &\sigma_0^2=1,
        &&\sigma_1^2=1, \\
      &[\sigma_0, \sigma_2]=[\sigma_0, \sigma_3]=1,
	&&[\sigma_1, \sigma_3]=[\sigma_1, \sigma_4]=1, \\
      &\sigma_0\sigma_1\sigma_0 = \sigma_1\sigma_0\sigma_1,
	&&\sigma_1\sigma_2\sigma_1 = \sigma_2\sigma_1\sigma_2, \\
      &\lambda_1^{-1}\sigma_2\lambda_1 = \sigma_3,
        &&\lambda_1^{-1}\sigma_3\lambda_1 = \sigma_4, \\
     &\sigma_0\lambda_0 = \lambda_1\sigma_0\sigma_1,
        &&\sigma_1\lambda_1 = \lambda_2\sigma_1\sigma_2, \\
      &[\sigma_0, \lambda_2]=[\sigma_0, \lambda_3]=1,
	&&[\sigma_1, \lambda_3]=[\sigma_1, \lambda_4]=1
  \rangle, \\ \\
\widehat{BV} = 
  \langle
    \lambda_0, \sigma_0, \sigma_1 \mid
      &\lambda_1^{-1}\lambda_2\lambda_1 = \lambda_3,
        &&\lambda_1^{-1}\lambda_3\lambda_1 = \lambda_4, \\
      &[\sigma_0, \sigma_2]=[\sigma_0, \sigma_3]=1,
	&&[\sigma_1, \sigma_3]=[\sigma_1, \sigma_4]=1, \\
      &\sigma_0\sigma_1\sigma_0 = \sigma_1\sigma_0\sigma_1,
	&&\sigma_1\sigma_2\sigma_1 = \sigma_2\sigma_1\sigma_2, \\
      &\lambda_1^{-1}\sigma_2\lambda_1 = \sigma_3,
        &&\lambda_1^{-1}\sigma_3\lambda_1 = \sigma_4, \\
      &\sigma_0\lambda_0 = 
          \lambda_1\sigma_0\sigma_1,
        &&\sigma_1\lambda_1 = 
          \lambda_2\sigma_1\sigma_2, \\
      &\sigma_0^{-1}\lambda_0 = 
          \lambda_1\sigma_0^{-1}\sigma_1^{-1},
        &&\sigma_1^{-1}\lambda_1 = 
          \lambda_2\sigma_1^{-1}\sigma_2^{-1}, \\
      &[\sigma_0, \lambda_2]=[\sigma_0, \lambda_3]=1,
	&&[\sigma_1, \lambda_3]=[\sigma_1, \lambda_4]=1
  \rangle.
\end{alignat*}
\end{lemma}

\section{Infinite presentations for \protect\(V\protect\) and
\protect\(BV\protect\)}

\subsection{Generators}  Generators are easy to find.

\newcommand{\opi}{\overline{\pi}}

\begin{lemma}\mylabel{FirstBVGens} The union of the three infinite
sets \[\begin{split} \{v_n
&=\lambda_0^{n+1}\lambda_1\lambda_0^{-n-2} \mid n\in \N\}, \\
\{\pi_n &= \lambda_0^{n+2}\sigma_1 \lambda_0^{-n-2} \mid n\in\N\},
\\ \{\opi_n &= \lambda_0^{n+1} \sigma_0 \lambda_0^{-n-1} \mid n\in
\N\} \end{split}\] generates \(V\) and \(BV\).  \end{lemma}

\begin{proof} Assume that we can show that any
\(u=W\beta\lambda_0^{-n}\) is a product of the claimed generators
when \(W\) has length \(n\) and satisfies \tref{SimpleCond} and
\(\beta\) is in \(S_{n+1}\) or \(B_{n+1}\).  Then any
\((F\beta)(G\gamma)^{-1}\) satsifying
\tref{BalancedCond}--\tref{BraidCond} can be written as \((F\beta
\lambda_0^{-n})(G\gamma\lambda_0^{-n})^{-1}\) and will be a product
of the claimed generators.  Thus we concentrate on a word such as
\(u\).

We first assume that \(\beta=1\) and argue inductively on the length
of \(W\) that \(W\lambda_0^{-n}\) is product of the generators in
\(\{v_0, v_1, \ldots\}\).  If \(W\) has length \(n\) and satisfies
\tref{SimpleCond}, then \(W\lambda_i\) satsfies \tref{SimpleCond} if
and only if \(0\le i\le n\).  Now if \(W\lambda_0^{-n}\) is a
product of the \(v_j\), then so is \(W\lambda_0\lambda_0^{-n-1}\),
so we may as well assume \(1\le i\le n\) which gives \(1-i\le0\).
This leads to the calculation \[\begin{split} W\lambda_0^{-n}v_{n-i}
&=W\lambda_0^{-n}\lambda_0^{n-i+1} \lambda_1 \lambda_0^{i-n-2} \\ &=
W\lambda_0^{1-i}\lambda_1 \lambda_0^{i-n-2} \\ &=W\lambda_i
\lambda_0^{-n-1} \end{split}\]  which finishes the induction.

For \(\beta\) non-trivial, we will have our result inductively if
\(W\beta\sigma_i \lambda_0^{-n} =ug\) for some generator \(g\)
whenever \(0\le i< n\).  Note that there is nothing to show unless
\(n\ge1\).  For \(i=0\), we have \(W\beta\lambda_0^{-n}\opi_{n-1} =
W\beta\sigma_0 \lambda_0^{-n}\).  For \(0<i<n\), we have
\[\begin{split} W\beta\lambda_0^{-n}\pi_{n-i-1} &=
W\beta\lambda_0^{-n} \lambda_0^{n-i+1}\sigma_1 \lambda_0^{-n+i-1} \\
&= W\beta\lambda_0^{1-i}\sigma_1 \lambda_0^{-n+i-1} \\ &=
W\beta\sigma_i\lambda_0^{-n}. \end{split} \] This completes the
proof.  \end{proof}

\subsection{Relations} The groups \(V\) and \(BV\) sit inside groups
with a known set of relations.  It is not hard to find relations for
the smaller groups.  The problem is finding enough relations among
them to get a presentation.  We start with lists of relations that
are demonstrably true.

\begin{lemma}\mylabel{InfiniteBVRels} In the following list, the
relations \tref{InfiniteBVRelsI} through \tref{InfiniteBVRelsX} hold
in \(BV\) and all the relations below hold in \(V\).
\mymargin{InfiniteBVRelsI--X \\ InfiniteVRelsI--II}
\begin{alignat}{2} 
v_qv_m&=v_mv_{q+1}, 
  &\qquad&m<q, \label{InfiniteBVRelsI}\\ 
\pi_q v_m&=v_m\pi_{q+1}, 
  && m<q, \label{InfiniteBVRelsII}\\ 
\pi^\epsilon_mv_m &= v_{m+1}\pi^\epsilon_m\pi^\epsilon_{m+1}, 
  &&m\ge0,\,\,\epsilon=\pm1, \label{InfiniteBVRelsIII}\\
\pi_qv_m &=v_m\pi_q, 
  &&m>q+1, \label{InfiniteBVRelsIV}\\
\opi_qv_m &=v_m\opi_{q+1}, 
  &&m<q, \label{InfiniteBVRelsV}\\
\opi^{\,\epsilon}_mv_m &= \pi^\epsilon_m\opi^{\,\epsilon}_{m+1}, 
  &&m\ge0,\,\,\epsilon=\pm1, \label{InfiniteBVRelsVI}\\
\pi_q\pi_m &= \pi_m\pi_q, 
  &&|m-q|\ge2, \label{InfiniteBVRelsVII}\\
\pi_m\pi_{m+1}\pi_m &= \pi_{m+1}\pi_m\pi_{m+1},
  &&m\ge0, \label{InfiniteBVRelsVIII}\\
\opi_q\pi_m &= \pi_m\opi_q, 
  &&q\ge m+2, \label{InfiniteBVRelsIX}\\
\pi_m\opi_{m+1}\pi_m &= \opi_{m+1}\pi_m\opi_{m+1},
  &&m\ge0,\label{InfiniteBVRelsX}  \\
\pi_m^2 &= 1,
  &&m\ge0,\label{InfiniteVRelsI} \\
\opi_m^2 &= 1,
  &&m\ge0.\label{InfiniteVRelsII}
\end{alignat}  
\end{lemma}

\begin{proof} It is straightforward to verify the claimed relations
from the definitions in Lemma \ref{FirstBVGens} and the relations
that hold in \(\widehat{V}\) and \(\widehat{BV}\) as given in
\tref{InfinitePresI} and \tref{InfinitePresII}.   We verify
\tref{InfiniteBVRelsI} as an 
example.   We assume \(m<q\) and  calculate 
\[\begin{split} v_qv_m &= \lambda_0^{q+1} \lambda_1 \lambda_0^{-q-2}
\lambda_0^{m+1} \lambda_1 \lambda_0^{-m-2} \\ &= \lambda_0^{q+1}
\lambda_1 \lambda_0^{m-q-1} \lambda_1 \lambda_0^{-m-2} \\ &=
\lambda_0^{q+1} \lambda_1 \lambda_{q-m+2} \lambda_0^{-q-3}, \\ \\
v_mv_{q+1} &= \lambda_0^{m+1} \lambda_1 \lambda_0^{-m-2}
\lambda_0^{q+2} \lambda_1 \lambda_0^{-q-3} \\ &= \lambda_0^{m+1}
\lambda_1 \lambda_0^{q-m} \lambda_1 \lambda_0^{-q-3} \\ &=
\lambda_0^{q+1} \lambda_{1+q-m} \lambda_1 \lambda_0^{-q-3} \\ &=
\lambda_0^{q+1} \lambda_1 \lambda_{q-m+2}
\lambda_0^{-q-3}. \end{split}\]
\end{proof}

The relations above are highly redundant.  In fact, the family of
generators is also overly large since the relations
\tref{InfiniteBVRelsVI} can be used to define each \(\pi_n\) as
\(\opi_n v_n \opi_{n+1}^{\,-1}\).  However, the large supply of
relations will make them easy to work with.

\subsection{A subgroup isomorphic to \protect\(F\protect\)}

Understanding words only in the \(v_i\) will be important.

\begin{lemma}\mylabel{FsubGp} The subgroup of \(V\) or \(BV\)
generated by \(\{v_i\mid i\in \N\}\) is isomorphic to the group
\(F\) with presentation \ref{GroupFPres}.  \end{lemma}

\begin{proof} By comparing \tref{InfiniteBVRelsI} with
\tref{GroupFPres}, we see that the subgroup in
question is a quotient of \(F\).  By \tref{NPNAQ}, \(F\) has no proper
non-abelian quotients, so we need only show \(v_2=v_0^{-1}v_1v_0\ne 
v_1\).  But if \(v_1=v_2\), then \[v_2v_1^{-1} =
(\lambda_0^3\lambda_1\lambda_0^{-4})
(\lambda_0^3\lambda_1^{-1}\lambda_0^{-2}) =
\lambda_0^3\lambda_1\lambda_2^{-1} \lambda_0^{-3}\] is trivial and
so we would have \(\lambda_1=\lambda_2\).  But the normal form given
in \tref{ForNorm} forbids this equality.
\end{proof}

\subsection{Calculations from the relations for
\protect\(BV\protect\)} 

From this point we derive consequences of the relations in Lemma
\ref{InfiniteBVRels}.  In order to insure that the work we do
applies equally to \(V\) and \(BV\), we will avoid making any use at
all of the relations \tref{InfiniteVRelsI} and
\tref{InfiniteVRelsII} until we loudly announce that we are ready to
resume using them.  Given two words \(w\) and \(w'\) in the
generators of Lemma \ref{FirstBVGens}, we will say \(w\sim w'\) to
mean that \(w\) can be converted to \(w'\) by use of the relations
\tref{InfiniteBVRelsI} through \tref{InfiniteBVRelsX} of Lemma
\ref{InfiniteBVRels}.

There is a function from words in the \(\pi_i\) and their inverses
into the group \(S_\infty\) of finitely supported permutations of
\(\N\).  In fact, it is a homomorphism from the subgroup of \(BV\)
generated by the \(\pi_i\), but we do not know this now, and have no
need to know it.  The function is gotten by taking \(\pi_i\) to the
transposition that interchanges \(i\) and \(i+1\).  Under this
assignment, any word \(w\) in the \(\pi_i\) and their inverses gives
a permutation that we continue to call \(w\).  This is used in the
next lemma.

\begin{lemma}\mylabel{PiVMuActs} If \(W\) is the set of words in the
\(\pi_i\) and their inverses, and if \(V\) is the set \(\{v_0, v_1,
\ldots\}\), then there are functions \((w,v)\mapsto w^v\) from
\(W\times V\rightarrow W\) and \((v,w)\mapsto v\cdot w\) from
\(V\times W\rightarrow W\) so that for \(w\in W\) and \(v_m\in V\),
we have \(wv_m \sim  v_j w^{v_m}\) where \(j=w(m)\) and \(v_m^{-1}w \sim 
(v_m\cdot w)v_k^{-1}\) where \(k=w^{-1}(m)\).

Further, if \(w\) is a word in \(\{\pi_i, \pi_i^{-1} \mid 0\le i\le
k\}\), then \(v_m\cdot w\) and \(w^{v_m}\) are words in \(\{\pi_i,
\pi_i^{-1} \mid 0\le i\le k+1\}\) and \(w(m)=m\) if \(m>k+1\) and
\(w(m)\le k+1\) if \(m \le k+1\).  \end{lemma}

\begin{proof}
From \tref{InfiniteBVRelsIII}, we get \(v_{m+1}^{-1}\pi^\epsilon_m
\sim \pi^\epsilon_m\pi^\epsilon_{m+1} v_m^{-1}\) which on inversion
gives \mymargin{InfiniteBVRelsXI} \begin{equation}
\label{InfiniteBVRelsXI} \pi^{-\epsilon}_mv_{m+1} \sim v_m
\pi^{-\epsilon}_{m+1}\pi^{-\epsilon}_m.\end{equation} This combines
with \tref{InfiniteBVRelsII}, \tref{InfiniteBVRelsIII} and
\tref{InfiniteBVRelsIV} and an induction argument to give \(wv_m
\sim v_j w^{v_m}\).  Inverting \tref{InfiniteBVRelsIII}, conjugating
\tref{InfiniteBVRelsII} and \tref{InfiniteBVRelsIV} by \(v_m\) and
combining with \(v_{m+1}^{-1}\pi^\epsilon_m \sim
\pi^\epsilon_m\pi^\epsilon_{m+1} v_m^{-1}\) and an induction
argument gives \(v_m^{-1}w \sim (v_m\cdot w)v_k^{-1}\).  The last
paragraph of the statement comes from the forms of
\tref{InfiniteBVRelsII}, \tref{InfiniteBVRelsIII},
\tref{InfiniteBVRelsIV} and \tref{InfiniteBVRelsXI} and the fact
that under the hypotheses of that paragraph, \(w\) maps to a
permutation in \(S_\infty\) that fixes all \(j\in \N\) with
\(j>k+1\).  \end{proof}

The interaction of \(v_m\) and \(\opi_q\) are partly covered by
\tref{InfiniteBVRelsV} and \tref{InfiniteBVRelsVI}.  The next lemma
covers the rest.

\begin{lemma}\mylabel{OpiVMuAct} The following hold for all \(k>0\):
\[\begin{split} \opi^{\,\epsilon}_m v_{m+k} &\sim  (v_mv_{m+1}\cdots
v_{m+k-2}v^2_{m+k-1}) \opi^{\,\epsilon}_{m+k+1} (\pi^\epsilon_{m+k}
\pi^\epsilon_{m+k-1} \cdots \pi^\epsilon_m), \\ v^{-1}_{m+k}
\opi^{\,\epsilon}_m &\sim  (\pi_m^\epsilon \pi_{m+1}^\epsilon \cdots
\pi^\epsilon_{m+k}) \opi^{\,\epsilon}_{m+k+1} (v_mv_{m+1}\cdots
v_{m+k-2}v^2_{m+k-1})^{-1}. \end{split}\] \end{lemma}

\begin{proof} From the inverted form of \tref{InfiniteBVRelsVI}, we
get \(\opi^{\,\epsilon}_m \sim  v_m \opi^{\,\epsilon}_{m+1} \pi_m^\epsilon \),
and inductively, get \[\opi^{\,\epsilon}_m \sim  (v_mv_{m+1}\cdots
v_{m+k-2}v_{m+k-1}) \opi^{\,\epsilon}_{m+k} (\pi^\epsilon_{m+k-1}
\pi^\epsilon_{m+k-2} \cdots \pi^\epsilon_m).\] Now \[\begin{split}
\opi^{\,\epsilon}_m &v_{m+k} \\ &\sim  (v_mv_{m+1}\cdots v_{m+k-2}v_{m+k-1})
\opi^{\,\epsilon}_{m+k} (\pi^\epsilon_{m+k-1} \pi^\epsilon_{m+k-2}
\cdots \pi^\epsilon_m) v_{m+k} \\ &\sim  (v_mv_{m+1}\cdots
v_{m+k-2}v_{m+k-1}) \opi^{\,\epsilon}_{m+k} \pi^\epsilon_{m+k-1}
v_{m+k}( \pi^\epsilon_{m+k-2} \pi^\epsilon_{m+k-3} \cdots
\pi^\epsilon_m) \\ &\sim  (v_mv_{m+1}\cdots v_{m+k-2}v_{m+k-1})
\opi^{\,\epsilon}_{m+k} v_{m+k-1} (\pi^\epsilon_{m+k}
\pi^\epsilon_{m+k-1} \cdots \pi^\epsilon_m) \\ &\sim  (v_mv_{m+1}\cdots
v_{m+k-2}v_{m+k-1}) v_{m+k-1} \opi^{\,\epsilon}_{m+k+1}
(\pi^\epsilon_{m+k} \pi^\epsilon_{m+k-1} \cdots \pi^\epsilon_m)
\end{split}\] where we have used \tref{InfiniteBVRelsIV},
\tref{InfiniteBVRelsXI} and \tref{InfiniteBVRelsV}.  The other
formula is a similar exercise starting with the uninverted form of
\tref{InfiniteBVRelsVI}.  \end{proof}

\begin{lemma}\mylabel{BVFirstForm} Let \(W\) be a word in the
generators of Lemma \ref{FirstBVGens}.  Then \(W\sim LMR\) where
\(L\) is a word in \(\{v_i\mid i\in \N\}\), \(R\) is a word in
\(\{v_i^{-1}\mid i\in \N\}\) and \(M\) is a word in \(\{\pi_i,
\pi_i^{-1} \opi_i, \opi_i^{\,-1} \mid i \in \N\}\).  \end{lemma}

\begin{proof} We start by assuming that \(x\) is represented by a
word in \(\{v_i, v_i^{-1}\mid i \in \N\}\).  From
\tref{InfiniteBVRelsI}, we get \(v_m^{-1} v_q^\epsilon \sim 
v_{q+1}^\epsilon v_m^{-1}\) and \(v_q^\epsilon v_m \sim  v_m
v_{q+1}^\epsilon\) whenever \(m<q\).  This allows negative powers
with low subscripts to pass to the right and positive powers with
low subscripts to pass to the left, each at the expense of raising
the already higher subscript.  When there is a possible ambiguity
(negative power wants to pass to the right of a positive power),
then one of the subscripts will be lower and the ambiguity is
resolved, or we have \(v_i^{-1}v_i\) which is replaced by 1.  Thus
we can achieve the promised form with \(M=1\).

Next we assume that \(x\) is represented by a word in positive and
negative powers of the generators of Lemma \ref{FirstBVGens} in
which there are no appearances of any \(\opi_i\) or
\(\opi_i^{\,-1}\).  Now Lemma \ref{PiVMuActs} allows us to pass
negative powers of the \(v_i\) from left to right and positive
powers of \(v_i\) from right to left over appearances of the
\(\pi_i^\epsilon\).  This comes at a cost of possibly increasing the
number of appearances of the \(\pi_i^\epsilon\) (if
\tref{InfiniteBVRelsIII} is used), but the number of
\(v_i^\epsilon\) cannot go up.  An easy induction on a complexity
derived from the word representing \(x\) completes the argument.
This is essentially the argument that would go to show that a
rewriting system derived from the relevant relations is terminating.

Now we assume that \(x\) is represented by an arbitrary word in the
positive and negative powers of the generators of Lemma
\ref{FirstBVGens}.  We will use Lemma \ref{OpiVMuAct} to pass
appearances of the \(v_i^\epsilon\) over appearances of the
\(\opi_j^{\,\delta}\).  Unfortunately, the number of appearances of
the \(v_i^\epsilon\) might go up, but the number of
\(\opi_j^{\,\delta}\) will not.

Consider the rightmost appearance of some \(\opi_j^{\,\delta}\) in
the word representing \(x\).  From the previous paragraphs, we can
assume that all appearances of positive powers of the \(v_i\) that
appear somewhere to the right of this \(\opi_j^{\,\delta}\) are
immediately to its right.  Now an application of
\tref{InfiniteBVRelsV}, \tref{InfiniteBVRelsVI} or the first
relation in Lemma \ref{OpiVMuAct} reduces by one the number of
positive powers of \(v_i\) to the right of this (altered) rightmost
appearance of \(\opi_{j+1}^{\,\delta}\) in the word.  If Lemma
\ref{OpiVMuAct} is used, then extra appearances of the \(\pi_k\) are
introduced, but these can be gotten ``out of the way'' by the
previous paragraphs.  Eventually, there will be no appearances of
positive powers of the \(v_i\) to the right of the rightmost
appearance of an \(\opi_n^{\,\delta}\) in the word.  An induction on
a complexity derived from the word shows that all positive powers of
the \(v_i\) can be moved to the extreme left of the word.  Now
negative powers can be moved to the right using the second relation
from Lemma \ref{OpiVMuAct} and altered forms of
\tref{InfiniteBVRelsV} and \tref{InfiniteBVRelsVI}.  \end{proof}

We now concentrate on the subword \(M\) obtained from Lemma
\ref{BVFirstForm}.  It is a word in \(\{\pi_i, \pi_i^{-1} \opi_i,
\opi_i^{\,-1} \mid i \in \N\}\).  We can map \(M\) to a word in the
generators of \(\widehat{BV}\) using the definitions in Lemma
\ref{FirstBVGens}.  Under the mapping, each letter in \(M\) becomes
some \(\sigma_i^\epsilon\) conjugated by a power of \(\lambda_0\).
We will be happiest when all the conjugations are by the same power
of \(\lambda_0\).  The appearances of \(\pi_i\) in \(M\) give us
more flexibility in achieving this because of the following
observation in which \(k\ge0\): \[\begin{split} \pi_i &=
\lambda_0^{i+2}\sigma_1 \lambda_0^{-i-2} \\ &= \lambda_0^{i+2}
\lambda_0^k \lambda_0^{-k} \sigma_1 \lambda_0^{-i-2} \\ &=
\lambda_0^{(i+k)+2} \sigma_{k+1} \lambda_0^{-(i+k)-2}. \end{split}\]
However, the \(\opi_i\) are more complicated.  This motivates the
next definitions.

Because of the definition \(\opi_n = \lambda_0^{n+1} \sigma_0
\lambda_0^{-n-1}\), we define the {\itshape height} of \(\opi_n\) to
be the subset \(\{n+1\}\sub \N\).  We define the {\itshape height}
of \(\pi_n\) to be the subset \(\{j\mid j\ge n+2\}\sub\N\).  If
\(M\) is a word in \(\{\pi_i, \pi_i^{-1} \opi_i, \opi_i^{\,-1} \mid
i \in \N\}\), then we define the height of \(M\) to be the
intersection of the heights of all the letters in \(M\).  Note that
it can easily be that the height of \(M\) is empty.  One of our
goals is the next lemma.

\begin{lemma}\mylabel{BVSecondForm} Let \(W\) be a word in the
generators of Lemma \ref{FirstBVGens}.  Then \(W\sim LMR\) as
specified in Lemma \ref{BVFirstForm} so that the height of \(M\) is
not empty.  \end{lemma}

We need some technical information about the behavior of height.
This will help in the proof and application of Lemma
\ref{BVSecondForm}.  In the following lemmas, a {\itshape
monosyllable} is a word in \(\{\pi_i, \pi_i^{-1} \opi_i,
\opi_i^{\,-1} \mid i \in \N\}\) with exactly one appearance of
\(\opi_i^{\,\pm1}\).

\begin{lemma} \mylabel{BVMonoHeight}Let \(M\) be a monosyllable of
height \(\{h\}\) and let \(0\le m<h\).  Then the following are true
where \(M'\) is some monosyllable of height \(\{h+1\}\) and \(0\le
j<h\).  The different appearances of \(M'\) and \(j\) are not to be
taken as representing the same values. {\claimenum \item \(M\sim
M'v_j^{-1}\), \item \(M\sim v_jM'\), \item \(M v_m \sim M'\) or \(M
v_m\sim v_j M'\), \item \(v_m^{-1} M\sim M'\) or \(v_m^{-1} M \sim
M' v_j^{-1}\).  \claimenumend}\end{lemma}

\begin{proof} The word \(M\) has the form
\(\Sigma_1\opi_{h-1}^{\,\gamma}\Sigma_2\) where
\(\gamma\in\{-1,1\}\) and \(\Sigma_1\) and \(\Sigma_2\) are words in
\(\{\pi_i, \pi_i^{-1} \mid 0\le i\le h-2\}\).

To prove (a), we get \(\opi_{h-1}^{\,\gamma} \sim \pi_{h-1}^\gamma
\opi_h^{\,\gamma} v_{h-1}^{-1}\) from \tref{InfiniteBVRelsVI}, and
from Lemma \ref{PiVMuActs}, we get \(v_{h-1}^{-1}\Sigma_2 \sim
\Sigma'_2v_j\) where \(0\le j <h\) and \(\Sigma'_2\) is a word in
\(\{\pi_i, \pi_i^{-1} \mid 0\le i\le h-1\}\).  Now \(M\sim \Sigma_1
\pi_{h-1}^\gamma \opi_h^{\,\gamma} \Sigma'_2v_j\).  The proof of (b)
is similar.

We consider (c).  By Lemma \ref{PiVMuActs}, \(\Sigma_2v_m\sim
v_k\Sigma'_2\) where \(0\le k<h\) and \(\Sigma'_2\) is a word in
\(\{\pi_i, \pi_i^{-1} \mid 0\le i\le h-1\}\).  If \(k=h-1\), then
\(\opi_{h-1}^{\,\gamma} v_{h-1} \sim
\pi^\gamma_{h-1}\opi_h^{\,\gamma}\) and we have \(Mv_m\sim \Sigma_1
\pi^\gamma_{h-1}\opi_h^{\,\gamma} \Sigma'_2\).  If \(k<h-1\), then
\(\opi_{h-1}^{\,\gamma} v_k \sim v_k \opi_h^{\,\gamma}\) and, by
Lemma \ref{PiVMuActs}, \(\Sigma_1 v_k \sim v_j \Sigma'_1\) where
\(0\le j<h\) and \(\Sigma'_1\) is a word in \(\{\pi_i, \pi_i^{-1}
\mid 0\le i\le h-1\}\).  Now \(Mv_m \sim
v_j\Sigma'_1\opi_n^{\,\gamma} \Sigma'_2\).  The proof of (b) is
similar.  \end{proof}

\begin{lemma}\mylabel{BVRaiseHeightI} Let \(W=M_1M_2\cdots M_t\)
where each \(M_i\) is a monosyllable of height \(\{h_i\}\) and
\(h_1\le h_2\le \cdots \le h_t\).  Then \(W\sim M'_1M'_2\cdots M'_t
v_j^{-1}\) or \(W\sim M'_1 M'_2 \cdots M'_t\) where \(0\le j<h_t\)
and each \(M'_i\) is a monosyllable of height \(\{h_i+1\}\).
\end{lemma}

\begin{proof} This is a repetitive application of Lemma
\ref{BVMonoHeight}.  We start by replacing \(M_1\) by
\(M'_1v_{i_1}^{-1}\) using Lemma \ref{BVMonoHeight}(a), and looking
at \(v_{i_1}^{-1}M_2\).  This either gives \(M'_2\) or
\(M'_2v_{i_2}^{-1}\) by Lemma \ref{BVMonoHeight}(d).  In the first
case, we continue by applying Lemma \ref{BVMonoHeight}(a) to
\(M_3\), and in the second case, we continue by applying Lemma
\ref{BVMonoHeight}(d) to \(v_{i_2}^{-1}M_3\).  Eventually, we
finish.  \end{proof}

\begin{lemma}\mylabel{BVRaiseHeightII} Let \(M\) be a word in
\(\{\pi_i, \pi_i^{-1} \opi_i, \opi_i^{\,-1} \mid i \in \N\}\) with a
non-empty height that contains \(h\in \N\).  Then there are \(i<h\),
\(j<h\), \(\epsilon\in\{0,1\}\) and \(\delta\in\{0,1\}\) so that
\(M\sim v_i^\epsilon M'\) and \(M\sim M''v_j^{-\delta}\) so that
\(M'\) and \(M''\) are words in \(\{\pi_i, \pi_i^{-1} \opi_i,
\opi_i^{\,-1} \mid i \in \N\}\) each with non-empty height
containing \(h+1\).  \end{lemma}

\begin{proof} If \(M\) contains no appearances of any
\(\opi_i^{\,\gamma}\), then the height of \(M\) is an infinite
subset of \(\N\) and we can let \(\epsilon\) and \(\delta\) be 0 and
let \(M'=M''=M\).

If there are appearances of the \(\opi_i^{\,\gamma}\) in \(M\), then
\(M\) is a concatenation of monosyllables of constant height.  We
get the form \(M''v_j^{\delta}\) by a direct application of Lemma
\ref{BVRaiseHeightI} and we get the form \(v_i^\epsilon M'\) by
applying Lemma \ref{BVRaiseHeightI} to \(M^{-1}\) and then inverting
the result.  \end{proof}

\begin{proof}[Proof of Lemma \ref{BVSecondForm}] Let \(W\sim LMR\)
as given by Lemma \ref{BVFirstForm}.  We may assume that \(M\) has
at least one appearance of some \(\opi_i^{\,\pm1}\).  Thus \(M\) is
a concatenation \(M_1M_2\cdots M_t\) of monosyllables.  Let
\(\{h_i\}\) be the height of the monosyllable \(M_i\).  If we alter
\(M\) by the use of the Lemmas \ref{BVMonoHeight} and
\ref{BVRaiseHeightI}, then we continue to refer to the monosyllables
as \(M_i\) and the heights as \(\{h_i\}\).  If we apply Lemma
\ref{BVRaiseHeightI} repeatedly to \(M_t\) we can raise the height
of \(M_t\) possibly at the expense of introducing extra terms into
the subword \(R\).  Thus we can assume that \(h_t\ge h_i\) for all
\(i<t\).  Now we apply Lemma \ref{BVRaiseHeightI} to \(M_{t-1}M_t\)
repeatedly so that we can assume \(h_{t-1}\ge h_i\) for \(i<t-1\)
and \(h_{t-1}\le h_t\).  Eventually, we get \(h_1\le h_2 \le \cdots
\le h_t\).

Now we invert the word \(LMR\) which still has the form given by
Lemma \ref{BVFirstForm} and \(M\) is again a concatenation of
monosyllables \(M_1M_2\cdots M_t\) with associated heights
satisfying \(h_1\ge h_2\ge \cdots \ge h_t\).  We repeat the above
process, but now we can obtain \(h_{t-1}=h_t\) after applying Lemma
\ref{BVRaiseHeightI} to \(M_t\), and we get \(h_{t-2}=h_{t-1}=h_t\)
after applying Lemma \ref{BVRaiseHeightI} to \(M_{t-1}M_t\), and so
forth.  Eventually, we get equal heights for all the monosyllables
which completes the proof.  \end{proof}

We now consider all of \(LMR\) as obtained from Lemma
\ref{BVFirstForm} and its mapping to a word in the generators of
\(\widehat{BV}\).  From Lemma \ref{BVSecondForm}, we can assume that
\(M\) maps to a word of the form \(\lambda_0^h w\lambda_0^{-h}\)
where \(w\) is a word in \(\{\sigma_i, \sigma_i^{-1}\mid i \in
\N\}\).  Since \(L\) is a word in \(\{v_i\mid i\in \N\}\), we can
map \(L\) to a word in the form
\mymargin{HeightKForm}\begin{equation} \label{HeightKForm}
\lambda_{i_1}\lambda_{i_2}\cdots \lambda_{i_k}
\lambda_0^{-k}.\end{equation} That one value of \(k\) correctly
appears in two places in \tref{HeightKForm} follows easily from the
definitions of the \(v_i\) and from the relations in 
\tref{InfinitePresII}.  We are not
concerned with the uniqueness of 
\(k\) and will say the {\itshape height} of \(L\) is no more than
\(k\) if \(L\) maps to a word in the form of \tref{HeightKForm}.
Since \(R^{-1}\) is a word in \(\{v_i\mid i\in \N\}\), we can say
the {\itshape height} of \(R\) is no more than \(k\) if the height
of \(R^{-1}\) is no more than \(k\).  Our next goal is the following
improvement on Lemmas \ref{BVFirstForm} and \ref{BVSecondForm}.

\begin{lemma}\mylabel{BVThirdForm} Let \(W\) be a word in the
generators of Lemma \ref{FirstBVGens}.  Then there is a \(k\) and
\(W\sim LMR\) as specified in Lemma \ref{BVSecondForm} so that the
height of \(M\) contains \(k\) and the heights of \(L\) and \(R\)
are no more than \(k\).  \end{lemma}

We need one more lemma about the behavior of height.

\begin{lemma}\mylabel{BVRaiseHeightIII} Let \(L\) be a word in
\(\{v_i\mid i\in \N\}\) of height no more than \(k\).  Then \(Lv_m\)
has height no more than \(k+1\) if \(m\le k-1\) and has height no
more than \(m+2\) if \(m\ge k-1\).  \end{lemma}

\begin{proof} We have \(L=\Lambda\lambda_0^{-k}\) where \(\Lambda\)
is a word in \(\{\lambda_i\mid i \in \N\}\).  Let \(m=k-1+j\) so
that \(v_m=\lambda_0^{k+j}\lambda_1\lambda_0^{-k-j-1}\).  Now \(Lv_m
= \Lambda\lambda_0^j\lambda_1\lambda_0^{-k-j-1}\).  If \(j\ge0\),
then \(Lv_m\) has height no more than \(k+j+1=m+2\) and if \(j\le
0\), then \(Lv_m = \Lambda \lambda_{j+1}\lambda_0^{-k-1}\) and has
height no more than \(k+1\).  \end{proof}

\begin{proof}[Proof of Lemma \ref{BVThirdForm}]  Let \(W\sim LMR\)
as given by Lemma \ref{BVSecondForm}.  Let \(L\) have height no more
than \(k_1\) and let \(R\) have height no more than \(k_2\).  If
\(M\) has no appearances of any \(\opi_i^{\,\pm1}\), then the height
of \(M\) is an infinite set and there is nothing to prove.  Thus
assume that \(M\) is a concatenation of monosyllables and has height
\(\{h\}\).

Assume first that \(k_1\le h\) and \(k_2>h\).  Use Lemma
\ref{BVRaiseHeightII} to replace \(M\) by \(v_m^\epsilon M'\) where
\(m<h\), where \(\epsilon\in\{0,1\}\) and where \(M'\) has height
\(\{h+1\}\).  If \(\epsilon=0\), then the height of
\(Lv_m^\epsilon\) is no more than \(k_1\) which is less than
\(h+1\).  Otherwise from Lemma \ref{BVRaiseHeightIII}, the height of
\(Lv_m^\epsilon\) is no more than \(k_1+1\) or \(m+2\), both of
which are no more than \(h+1\).  Thus repeated uses of Lemma
\ref{BVRaiseHeightII} give us the desired result.  If \(k_1>h\) and
\(k_2\le h\), then we use Lemma \ref{BVRaiseHeightII} repeatedly
replacing \(M\) by \(M'v_m^{-\delta}\) with \(\delta\in\{0,1\}\).

Now assume \(k_1>h\) and \(k_2>h\).  Now repeated uses of Lemma
\ref{BVRaiseHeightII} replacing \(M\) by \(v_m^\epsilon M'\) will
achieve \(h\ge k_2\) at the expense of raising \(k_1\).  Now we are
in one of the cases of the previous paragraph.  \end{proof}

\subsection{Infinite presentations}

We now are ready to resume using the relations \tref{InfiniteVRelsI} and
\tref{InfiniteVRelsII}  when dealing with \(V\).  In the following
\(\sim\) will use only relations \tref{InfiniteBVRelsI} through
\tref{InfiniteBVRelsX} in the case of \(BV\) and will use
\tref{InfiniteVRelsI} and \tref{InfiniteVRelsII} in addition in the
case of \(V\).  The difference in the application will be clear.

\begin{prop}\mylabel{FinalReduction} Let \(W\) be a word in the
geneators of Lemma \ref{FirstBVGens}.  If \(W\) represents the
trivial element of either \(V\) or \(BV\), then \(W\sim\phi\) where
\(\phi\) is the empty word.
\end{prop}

\begin{proof} We may assume that \(W\) has the form \(LMR\) as given
by Lemma \ref{BVThirdForm}.

We first argue that \(W\sim \phi\) follows if we can show \(M\sim
\phi\).  If \(M\sim \phi\), then \(W\sim LR\) which is a word in
\(\{v_i, v_i^{-1}\mid i\in \N\}\).  By Lemma \ref{FsubGp}, the
subgroup generated by \(\{v_i, v_i^{-1}\mid i\in \N\}\) is
isomorphic to the group \(F\) whose relations are precisely the
relations in \tref{InfiniteBVRelsI}.  Thus if \(LR\) is the the
identity, then \(W\sim\phi\).

It remains to show that \(M\sim \phi\).

We know that \(M\) has a height that includes some \(h\) and that
\(L\) and \(R\) have height no more than \(h\).  Thus the mapping of
\(W\) into the generators of \(\widehat{V}\) or  \(\widehat{BV}\)
takes the form 
\(\Lambda_1 \lambda_0^{-s}\lambda_0^h w\lambda_0^{-h} \lambda_0^t
\Lambda_2\) where \(\Lambda_1\) and \(\Lambda_2^{-1}\) are words in
\(\{\lambda_i\mid i \in \N\}\), \(w\) is a word in \(\{\sigma_i\mid
i\in\N\}\) and \(s\le h\) and \(t\le h\) both hold.  Thus the image
of \(W\) is really of the form \(\Lambda'_1 w \Lambda'_2\) where
\(w\) is the same as above and \(\Lambda'_1\) and
\((\Lambda'_2)^{-1}\) are words in \(\{\lambda_i\mid i\in\N\}\).

We can say more about the letters in \(w\).  The mapping of the
elements of \(M\) into \(w\) took each letter of \(M\) to some
\(\sigma_i\) conjugated by \(\lambda_0^h\).  In particular, each
\(\opi_i^{\,\pm1}\) in \(M\) is of the form
\(\opi_{h-1}^{\,\epsilon}\) with \(\epsilon\in\{-1,1\}\).  This maps to
\(\lambda_0^h\sigma_0^\epsilon\lambda_0^{-h}\).  Each
\(\pi_i^\epsilon\) must have \(i\le h-2\) and maps to
\(\lambda_0^{i+2}\sigma_1^\epsilon\lambda_0^{-i-2}\).  However as
noted before Lemma \ref{BVSecondForm}, it can also be mapped to
\(\lambda_0^{(i+k)+2}\sigma_{1+k}^\epsilon\lambda_0^{-(i+k)-2}\).
Setting \((i+k)+2=h\) gives \(1+k=(h-1)-i\).

The previous paragraph gives a simple formula for getting the word
\(w\) from \(M\).  We have that
\mymargin{TheXVer}\begin{equation}\label{TheXVer}
M=X_{i_1}^{\epsilon_1}X_{i_2}^{\epsilon_2}\cdots
X_{i_r}^{\epsilon_r}\end{equation} where each \(X\) is either the
symbol \(\opi\) or the symbol \(\pi\).  The resulting word \(w\) is
\mymargin{TheSigmaVer}\begin{equation}\label{TheSigmaVer}
\sigma_{j_1}^{\epsilon_1} \sigma_{j_2}^{\epsilon_2} \cdots
\sigma_{j_r}^{\epsilon_r} \end{equation} where each
\(j_k=(h-1)-i_k\).

We first argue that it suffices to show that \tref{TheSigmaVer}
reduces to the trivial word modulo the relations in \(\widehat{V}\)
or \(\widehat{BV}\) that relate only to the \(\sigma_i\):
specifically the relations \tref{PermRelsB} and \tref{PermRelsC} in
either case and \tref{PermRelsA}, in addition, in the case of
\(\widehat{V}\).

Assume a sequence of applications of \tref{PermRelsB},
\tref{PermRelsC} and perhaps \tref{PermRelsA} reduce
\tref{TheSigmaVer} to the empty word.  The relation \tref{PermRelsB}
reads \(\sigma_m\sigma_n=\sigma_n\sigma_m\) if \(|m-n|\ge2\).
However, \(\sigma_m\) comes from \(X_p\) and \(\sigma_n\) comes from
\(X_q\) and with \(p=(h-1)-m\) and \(q=(h-1)-n\).  But
\(|p-q|=|m-n|\ge2\) and the relations of Lemma \ref{InfiniteBVRels}
imply that \(X_pX_q=X_qX_p\) so the application of \tref{PermRelsB}
to \tref{TheSigmaVer} can be mirrored by a corresponding application
of either \tref{InfiniteBVRelsVII} or \tref{InfiniteBVRelsIX} to
\tref{TheXVer}.  After the corresponding applications, the new
version of \tref{TheXVer} maps to the new version of
\tref{TheSigmaVer} under the mapping given by the paragraphs above.
An identical argument says that an application of \tref{PermRelsC}
to \tref{TheSigmaVer} can be mirrored by a corresponding application
of either \tref{InfiniteBVRelsVIII} or \tref{InfiniteBVRelsX} to
\tref{TheXVer}.  Finally, applications of \tref{PermRelsA} can be
mirrored by applications of either \tref{InfiniteVRelsI} or
\tref{InfiniteVRelsII}.  Note that applications of
\tref{InfiniteVRelsI} or \tref{InfiniteVRelsII} will only be needed
if \tref{PermRelsA} is needed.  However, \tref{PermRelsA} only
applies to \(\widehat{V}\) and so \tref{InfiniteVRelsI} or
\tref{InfiniteVRelsII} will only be needed if we are dealing with
\({V}\).

Thus the reduction of \tref{TheSigmaVer} to
the empty word gives instructions to achieve a reduction of
\tref{TheXVer} to the trivial word that uses appropriate relations
in each case.

Now we must argue that \tref{TheSigmaVer} can be trivialized by
applications of \tref{PermRelsB}, \tref{PermRelsC} and perhaps
\tref{PermRelsA}.

We are assuming the word \(\Lambda'_1 w\Lambda'_2\) represents the
identity in \(\widehat{V}\) or \(\widehat{BV}\), that \(w\) is the
word \tref{TheSigmaVer} above, and that \(\Lambda'_1\) and
\((\Lambda'_2)^{-1}\) are words in \(\{\lambda_i\mid i \in\N\}\).
Thus \(\Lambda'_1w = (\Lambda'_2)^{-1}\) as elements of
\(\widehat{V}\) or \(\widehat{BV}\).  Since \(\widehat{V}\) and
\(\widehat{BV}\) are the groups of fractions of the monoid \(\FP\)
and \(\FB\), respectively, the embedding property \tref{GOFEmbed} of
groups of fractions says that \(\Lambda'_1w = (\Lambda'_2)^{-1}\) as
elements of \(\FP\) or \(\FB\).  From the uniqueness of
representation in \ZS products \tref{ZSUnique}, this says that
\(w=1\) as an element of \(S_\infty\) or \(B_\infty\).  Since the
relations \tref{PermRelsB} and \tref{PermRelsC} suffice to present
\(B_\infty\) and the relations \tref{PermRelsB}, \tref{PermRelsC}
and \tref{PermRelsA} suffice to present \(S_\infty\), we are done
with the proof.
\end{proof}

We are now ready to give an infinite presentation of \(V\) and \(BV\).

\begin{thm}\mylabel{InfiniteBVPres} The groups \(V\) and \(BV\) are
presented by by the generators of Lemma \ref{FirstBVGens} and the
relations \tref{InfiniteBVRelsI} through \tref{InfiniteBVRelsX} of
Lemma \ref{InfiniteBVRels} in the case of \(BV\) and the relations
\tref{InfiniteBVRelsI} through \tref{InfiniteVRelsII} of Lemma
\ref{InfiniteBVRels} in the case of \(V\).  \end{thm}

\begin{cora} \mylabel{InfiniteBVPresB}
The group \(BV\) is presented by the generators \(v_n\)
and \(\opi_n\) for \(n\in \N\) and by the relations 
\begin{alignat*}{2} 
v_qv_m&=v_mv_{q+1}, 
  &\qquad&m<q,\\
\pi^\epsilon_mv_m &= v_{m+1}\pi^\epsilon_m\pi^\epsilon_{m+1}, 
  &&m\ge0,\,\,\epsilon=\pm1,\\
\pi_qv_m &=v_m\pi_q, 
  &&m>q+1,\\
\opi_qv_m &=v_m\opi_{q+1}, 
  &&m<q, \\
\pi_m &= \opi^{\,-1}_{m+1}v_m^{-1}\opi_m, 
  &&m\ge0, \\
\pi_q\pi_m &= \pi_m\pi_q, 
  &&|m-q|\ge2, \\
\pi_m\pi_{m+1}\pi_m &= \pi_{m+1}\pi_m\pi_{m+1},
  &&m\ge0, \\
\opi_q\pi_m &= \pi_m\opi_q, 
  &&q\ge m+2, \\
\pi_m\opi_{m+1}\pi_m &= \opi_{m+1}\pi_m\opi_{m+1},
  &&m\ge0
\end{alignat*} 
in which \(\pi_n=\opi_nv_n\opi_{n+1}^{\,-1}\) is taken as a
definition for each \(n\in\N\).   

The group \(V\) is presented by the generators \(v_n\)
and \(\opi_n\) for \(n\in \N\) and by the relations with definitions
listed above for \(BV\) with the addition of the relations 
\[
\pi_m^2 = \opi_m^2=1,
  \qquad m\ge0.
\]
\end{cora}

\begin{proof} The definitions for the \(\pi_n\) are obtained from
the relations \tref{InfiniteBVRelsVI} with \(\epsilon=1\).  Now the
relations \tref{InfiniteBVRelsII} follow from the definitions of the
\(\pi_n\) and from \tref{InfiniteBVRelsI} and
\tref{InfiniteBVRelsV}.  The fifth line of relations above is
obtained from \tref{InfiniteBVRelsVI} with \(\epsilon=-1\).
\end{proof}

It is not clear if there is any way to conclude the relations
\(\pi_m^2=1\) from the others in the case of \(V\).

\subsection{A ``Coxeter group/Artin group'' pair}

The choice of generators for the presentations in Corollary
\ref{InfiniteBVPresB} was made mostly to mirror the traditional
presentations of Thompson's groups.  Instead, we could have turned
the definition \(\pi_n=\opi_nv_n\opi_{n+1}^{\,-1}\) into
\(v_n=\opi_n^{\,-1}\pi_n\opi_{n+1}\) to get the following.

\begin{cora} \mylabel{InfiniteBVPresC}
The group \(BV\) is presented by the generators \(\pi_n\)
and \(\opi_n\) for \(n\in \N\) and by the relations 
in Corollary
\ref{InfiniteBVPresB}
in which \(v_n=\opi_n^{\,-1}\pi_n\opi_{n+1}\) is taken as a
definition for each \(n\in\N\).   

The group \(V\) is presented by the generators \(\pi_n\)
and \(\opi_n\) for \(n\in \N\) and by the relations with definitions
listed above for \(BV\) with the addition of the relations 
\[
\pi_m^2 = \opi_m^2=1,
  \qquad m\ge0.
\]
\end{cora}

\section{Finite presentations for \protect\(V\protect\) and
\protect\(BV\protect\)}

We can now extract finite presentations from Corollary
\ref{InfiniteBVPresB}. 

\begin{thm}\mylabel{FiniteBVPres} The group \(BV\) has the
following presentation.  The equalities  \(v_i= v_0^{1-i}v_1v_0^{i-1}\)
and \(\opi_i= v_0^{1-i}\opi_1
v_0^{i-1}\) for all \(i\ge2\), and 
\(\pi_i=\opi_iv_i\opi_{i+1}^{\,-1}\) for all \(i\in\N\)
are taken as definitions.
\begin{alignat*}{2}
BV = 
  \langle
    v_0, v_1, \opi_0, \opi_1 \mid
      &v_2v_1 = v_1v_3,
	&\qquad&v_3v_1 = v_1v_4, \\
      &\opi_2v_1 = v_1\opi_3, &&\opi_3v_1 = v_1\opi_4, \\
      &\pi_0v_0 = v_1\pi_0\pi_1, &&\pi_1v_1 = v_2\pi_1\pi_2, \\
      &\pi_0^{-1}v_0 = v_1\pi_0^{-1}\pi_1^{-1},
        &&\pi_1^{-1}v_1 = v_2\pi_1^{-1}\pi_2^{-1}, \\
      &\pi_0v_2 = v_2\pi_0, &&\pi_0v_3 = v_3\pi_0, \\
      &\pi_1v_3 = v_3\pi_1, &&\pi_1v_4 = v_4\pi_1, \\
      &\pi_0 = \opi^{\,-1}_{1}v_0^{-1}\opi_0,
        &&\pi_1 = \opi^{\,-1}_{2}v_1^{-1}\opi_1, \\
      &\pi_0\pi_2 = \pi_2\pi_0, &&\pi_0\pi_3 = \pi_3\pi_0, \\
      &\pi_1\pi_3 = \pi_3\pi_1, &&\pi_1\pi_4 = \pi_4\pi_1, \\
      &\pi_0\pi_1\pi_0 = \pi_1\pi_0\pi_1, 
	&&\pi_1\pi_2\pi_1 = \pi_2\pi_1\pi_2, \\
      &\opi_2\pi_0 = \pi_0\opi_2, &&\opi_3\pi_0 = \pi_0\opi_3, \\
      &\opi_3\pi_1 = \pi_1\opi_3, &&\opi_4\pi_1 = \pi_1\opi_4, \\
      &\pi_0\opi_1\pi_0 = \opi_1\pi_0\opi_1, 
	&&\pi_1\opi_2\pi_1 = \opi_2\pi_1\opi_2
  \rangle.
\end{alignat*}

A presentation for the group \(V\) is obtained from that for \(BV\)
by adding the relations \(\opi_0^2=\opi_1^2=\pi_0^2=\pi_1^2=1\).

The generating set for the presentations of both \(V\) and \(BV\)
can be replaced by \(\{\pi_0, \pi_1, \opi_0, \opi_1\}\) if the
definitions above are replaced by the equalities \(\pi_i=
v_0^{1-i}\pi_1v_0^{i-1}\) and \(\opi_i= v_0^{1-i}\opi_1 v_0^{i-1}\)
for all \(i\ge2\), and \(v_i=\opi_i^{\,-1}\pi_i\opi_{i+1}\) for all
\(i\in\N\).  \end{thm}

\begin{proof} The first presentations for \(BV\) and \(V\) are
extracted from the presentations of Corollary \ref{InfiniteBVPresB}
using the arguments from Section \ref{ThompsonArgsSec} in the same
way that the presentations of Lemma \ref{FinitePres} are extracted
from those of \tref{InfinitePresI} and \tref{InfinitePresII}.  Each
family of relations in Corollary \ref{InfiniteBVPresB} gives rise to
2 or four relations in the above list.  The presentations of the
last paragraph can be shown equivalent via Tietze transformations to
the first presentations of the statement.  \end{proof}

Note that several relations from Theorem \ref{FiniteBVPres} become
redundant in the case of \(V\) since the \(\pi_m\) and \(\opi_m\)
become their own inverses.

\providecommand{\bysame}{\leavevmode\hbox to3em{\hrulefill}\thinspace}

\noindent Department of Mathematical Sciences

\noindent State University of New York at Binghamton

\noindent Binghamton, NY 13902-6000

\noindent USA

\noindent email: matt@math.binghamton.edu

\end{document}